\documentclass[12pt]{article}
\usepackage{latexsym,amssymb,amsmath} 
\newtheorem{theorem}{Theorem}[section]

\newtheorem{dfn}[theorem]{Definition}

\newtheorem{lem}[theorem]{Lemma}
\newtheorem{cor}[theorem]{Corollary}
\newtheorem{thm}[theorem]{Theorem}

\newcommand{\prf}{{\bf Proof}. }

 \newcommand{\minpunkt}{
\begin{picture}(12,6)(-4,-2)
  \put(-1, 0){ \makebox(0,0){$-$} }
  \put(-1, 2){ \makebox(0,0){$\cdot$} }
  \end{picture}}

 \def\provx#1#2#3#4{
\setbox1=\hbox{\kern1.5pt$\scriptstyle#3$}
\def\zeichen{#2}
\ifx\zeichen\empty\setbox0=\hbox to .75em{}\else\setbox0=\hbox
{\kern1.5pt$\scriptstyle#2$}\fi
\dimen1=\dp0 \ifdim \dimen1=0pt
\advance \dimen1 by 1.5ex \else \advance \dimen1 by 1.2ex
\fi\dimen3=2ex\dimen4=.5ex\ifdim \wd0<\wd1 \dimen2=\wd1 \else \dimen2=\wd0
\fi\hbox{$#1\hskip 5pt minus5pt\vrule height\dimen3
depth\dimen4\raise\dimen1\copy0\hskip-1\wd0 \lower\ht1
\copy1\hskip-1\wd1\vrule width\dimen2 height.7ex depth-.6ex\hskip3pt
minus1.5pt#4\hskip2pt plus2pt minus2pt$}}

\def\prov#1#2#3{
\setbox1=\hbox{\kern1.5pt$\scriptstyle#2$}
\def\zeichen{#1}
\ifx\zeichen\empty\setbox0=\hbox to .75em{}\else\setbox0=\hbox
{\kern1.5pt$\scriptstyle#1$}\fi
\dimen1=\dp0
\ifdim \dimen1=0pt
\advance \dimen1 by 1.5ex \else \advance \dimen1 by 1.2ex
\fi\dimen3=2ex\dimen4=.5ex\ifdim \wd0<\wd1 \dimen2=\wd1 \else \dimen2=\wd0
\fi\hbox{\hskip0pt plus 4pt
$\vrule height\dimen3
depth\dimen4\raise\dimen1\copy0\hskip-1\wd0
\lower\ht1\copy1\hskip-1\wd1\vrule width\dimen2 height.7ex depth-.6ex
\hskip3pt minus1.5pt#3\hskip2pt plus2pt minus2pt$}}

\def\prv#1#2{
\setbox1=\hbox{\kern1.5pt$\scriptstyle#2$}
\ifx\zeichen\empty\setbox0=\hbox to .75em{}\else\setbox0=\hbox
{\kern1.5pt$\scriptstyle#1$}\fi
\dimen1=\dp0 \ifdim \dimen1=0pt
\advance \dimen1 by 1.5ex \else \advance \dimen1 by 1.2ex
\fi\dimen3=2ex\dimen4=.5ex\ifdim \wd0<\wd1 \dimen2=\wd1 \else \dimen2=\wd0
\fi\hbox{\hskip.5em$\vrule height\dimen3
depth\dimen4\raise\dimen1\copy0\hskip-1\wd0
\lower\ht1\copy1\hskip-1\wd1\vrule width\dimen2 height.7ex depth-.6ex
\hskip3pt minus1.5pt$}}

\mathchardef\str='1066
\def\negprov#1#2#3{
\setbox1=\hbox{\kern1.5pt$\scriptstyle#2$}
\setbox4=\hbox{$\str$}
\def\zeichen{#1}
\ifx\zeichen\empty\setbox0=\hbox to 1em{}\else\setbox0=\hbox
{\kern1.5pt$\scriptstyle#1$}\fi
\dimen1=\dp0
\ifdim \dimen1=0pt
\advance \dimen1 by 1.5ex \else \advance \dimen1 by 1.2ex
\fi\dimen3=2ex\dimen4=.5ex\ifdim \wd0<\wd1 \dimen2=\wd1 \else \dimen2=\wd0
\fi\hbox{\hskip.5em$\kern-1.9pt\raise1pt\copy4\kern-\wd4\kern1.9pt\vrule height\dimen3
depth\dimen4\raise\dimen1\copy0\hskip-1\wd0
\lower\ht1\copy1\hskip-1\wd1\vrule width\dimen2 height.7ex depth-.6ex
\hskip3pt minus1.5pt#3\hskip2pt plus2pt minus2pt$}}

\def\goed#1{\setbox5=\hbox{$#1$}\dimen1=.25em \dimen2=\dimen1 \advance \dimen2
by -1pt\hbox{\raise.65\ht5 \hbox{\vrule height.5\ht5 depth0pt width.4pt\vrule
height.5\ht5 width\dimen1 depth-.48\ht5}\kern-\dimen2\copy5\kern-\dimen2
\raise.65\ht5 \hbox{\vrule height .5\ht5 width\dimen1 depth-.48\ht5\vrule
height.5\ht5 depth 0pt width.4pt}\hskip4pt plus2pt minus2pt}}

\def\mod#1#2{
\def\zeichen{#1}
\hbox{\hskip 2pt plus3pt minus 2pt\vrule width.5pt height2ex depth.5ex
\vbox{\ifx\zeichen\empty\hbox to .75em{}\else
\hbox{\kern1.5pt $\scriptstyle#1$}\fi
\kern2pt
\hrule
\kern1.7pt
\hrule\kern1.7pt}
\hskip3pt minus 2pt$#2$}\hskip2pt
plus3pt minus2pt}

\def\notmod#1#2{\hbox{\hskip 2pt plus 3pt minus 3pt\vrule width.5pt
height2ex depth.5ex
\vbox{\hbox{\kern1.5pt $\scriptstyle#1$}\kern3pt
\setbox0=\hbox{\kern2pt$\scriptstyle/$}
\hrule
\kern-1.7pt
\copy0
\kern-\ht0
\kern 1.7pt
\hrule\kern1.7pt}n
\hskip3pt minus 2pt$#2$}\hskip2pt
plus3pt minus2pt}
\def\sq{\hbox{\rlap{$\sqcap$}$\sqcup$}}
\def\qed{\ifmmode\sq\else{\unskip\nobreak\hfil\penalty50\hskip1em\null
\nobreak\hfil\sq\parfillskip=0pt\finalhyphendemerits=0\endgraf}\fi\medskip}

\def\lleq{\hbox{\hskip3pt minus3pt\kern1pt\lower4pt
\vbox{\hbox{$\scriptstyle\ll$}
\kern-7pt\hbox{\kern1pt$\scriptstyle=$}}\hskip3pt minus 3pt}}

\mathchardef\res='1152
\mathchardef\qin='1062
\mathchardef\qprec='1036
\mathchardef\qless='474
\mathchardef\dpkt='72

\def\A{{\rm\kern.22em
   \vrule width.02em
       height0.5ex depth 0ex
   \kern-.24em A}}
\def\B{{\rm I\kern-.25em B}}
\def\C{{\rm\kern.24em
   \vrule width.02em
       height1.4ex depth-.05ex
   \kern-.26em C}}
\def\D{{\rm I\kern-.25em D}}
\def\E{{\rm I\kern-.25em E}}
\def\F{{\rm I\kern-.25em F}}
\def\G{{\rm\kern.24em
   \vrule width.02em
       height1.4ex depth-.05ex
   \kern-.26em G}}
\def\H{{\rm I\kern-.25em H}}
\def\I{{\rm I\kern-.25em I}}
\def\J{{\rm\kern.19em
   \vrule width.02em
       height1.47ex depth 0ex
   \kern-.21em J}}
\def\K{{\rm I\kern-.25em K}}
\def\L{{\rm I\kern-.25em L}}
\def\M{{\rm I\kern-.23em M}}
\def\N{{\rm I\kern-.23em N}}
\def\O{{\rm\kern.24em
   \vrule width.02em
       height1.4ex depth-.05ex
   \kern-.26em O}}
\def\P{{\rm I\kern-.25em P}}
\def\Q{{\rm\kern.24em
   \vrule width.02em
       height1.4ex depth-.05ex
   \kern-.26em Q}}
\def\R{{\rm I\kern-.25em R}}
\def\S{{\rm\kern.18em
   \vrule width.02em
       height1.47ex depth-.9ex
\kern.12em
\vrule width.02em
    height0.7ex depth 0ex
\kern-.34em S}}
\def\T{{\rm\kern.45em
   \vrule width.02em
       height1.47ex depth 0ex
   \kern-.47em T}}
\def\U{{\rm\kern.30em
   \vrule width.02em
       height1.47ex depth-.05ex
   \kern-.32em U}}
\def\V{{\rm\kern.27em
   \vrule width.02em
       height1.47ex depth-.8ex
   \kern-.29em V}}
\def\W{{\rm\kern.25em
   \vrule width.02em
       height1.47ex depth-.9ex
   \kern.34em
   \vrule width.02em
       height1.47ex depth-.9ex
   \kern-.63em W}}
\def\X{{\rm\kern.30em
   \vrule width.02em
       height1.47ex depth-1ex
   \kern.12em
   \vrule width.02em
       height0.4ex depth 0ex
   \kern-.46em X}}
\def\Y{{\rm\kern.25em
   \vrule width.02em
       height1.0ex depth 0ex
   \kern-.27em Y}}
\def\Z{{\rm\kern.26em
   \vrule width.02em
       height0.5ex depth 0ex
   \kern.04em
   \vrule width.02em
       height1.47ex depth-1ex
   \kern-.34em Z}}

\newcommand{\ko}{\mathrm C}
\newcommand{\nf}{=_{NF}}
\newcommand{\Tu}[2]{ \hat{\mathrm T}^{\omega}_{#1}(#2)}
\newcommand{\Tus}[2]{\mathrm T^{#1}_{\omega}(#2)}
\newcommand{\Su}[3]{\mathrm S^{#1}_{#2}(#3)}
\newcommand{\PA}{{\mathbf{PA}}}
\newcommand{\PRA}{{\mathbf{PRA}}}

\newcommand{\ord}{\mathrm{ord}}
\newcommand{\RR}{{\mathcal{R}}}
\newcommand{\PRWO}{\mathrm{PRWO}}
\begin{document}
\title{Goodstein's theorem revisited}
\author{Michael Rathjen\\
 {\small {\it School of Mathematics,
  University of Leeds}} \\
{\small {\it Leeds, LS2 JT, England}}}
\date{}
\maketitle

\begin{abstract} Inspired by Gentzen's 1936 consistency proof, Goodstein found
a close fit between descending sequences of ordinals $<\varepsilon_0$ and sequences
of integers, now known as Goodstein sequences.
This article revisits Goodstein's 1944 paper.  In light of new historical
details found in a correspondence between Bernays and Goodstein,
we address the question of how close Goodstein came to proving an independence result for $\PA$. We also present an elementary proof of the fact that already the termination of all special Goodstein
sequences, i.e. those induced by the shift function,  is not provable in $\PA$.
This was first proved by Kirby and Paris in 1982, using techniques from the model theory of arithmetic.
The proof presented here arguably only uses tools that would have been available in the 1940's or 1950's. Thus we ponder the question whether striking independence results
could have been proved much earlier? In the same vein we also wonder whether
the search for strictly mathematical examples of an incompleteness in $\PA$ really
 attained its ``holy grail" status before the late 1970's.
Almost no direct moral is ever given; rather, the paper strives to lay out evidence for the reader to consider and have the reader form their own conclusions.
However, in relation to independence results, we think that both Goodstein and Gentzen are deserving of more credit.
 \end{abstract}

\section{History}
This paper grew out of a Goodstein lecture that I gave at the Logic Colloquium 2012 in Manchester. The lecture touched on many of Goodstein's papers, though, this article will just be concerned with his best known result \cite{goodstein} from 1944. Whilst reading \cite{goodstein}, I formed the overwhelming
impression that Goodstein came very close to proving an independence result.
Recent
archival studies by Jan von Plato (cf. \cite{plato}) which have brought to light a remarkable
correspondence between Goodstein and Bernays have now confirmed that
this impression was not unfounded.
Goodstein's paper \cite{goodstein} originally bore the title ``A note on Gentzen's theorem", thereby referring to Gentzen's 1936 paper \cite{gentzen} which proved the consistency of first order number theory\footnote{First order number theory or {\em reine Zahlentheorie} as it was called by Gentzen is essentially the same system as what is nowadays called Peano arithmetic, $\PA$.} by transfinite induction up to the ordinal $\varepsilon_0$.
He sent it to Church  in 1942 for publication in the JSL. Church sent the paper to Bernays for refereeing.
 Bernays then contacted
Goodstein directly and included a long list of remarks and suggestions in his letter \cite{bernays} dated
1 September 1942.
As a result of these comments, Goodstein altered his paper considerably and also changed the title to ``On the restricted ordinal theorem".
By the latter he meant the proposition that every strictly descending sequence of ordinals below $\varepsilon_0$ is necessarily finite. As Gentzen \cite{gentzen} showed,  this implies the consistency of first order number theory.
 Crucially in his paper Goodstein proved that this statement is equivalent to a statement $P$ about integers, now known as the termination of Goodstein sequences.
 From Bernays' letter it is clear that the original version of Goodstein's paper contained a claim about the
 unprovability of $P$ in number theory. Bernays in his letter correctly pointed out that $P$, being of $\Pi^1_1$ form, is not a statement that
 can be formalized in Gentzen's system of first order number theory  as it  talks about all descending sequences.
 \begin{quote} The system $\underline{A}$ cannot be exactly the system denoted by Gentzen as ``reine Zahlentheorie", since this one contains no function variables and so your theorem $P$ is not expressible in it. However the Gentzen proof surely can be extended to the case that free function variables are added to the considered formal system.
 \cite{bernays}\end{quote}
 Unfortunately Goodstein then removed the passage about the unprovability of $P$. He could have easily\footnote{This expression may be too strong since
 it assumes that Goodstein had penetrated
  the details of Gentzen's rather difficult paper \cite{gentzen}.}
 come up with an independence result for $\PA$ as Gentzen's proof only utilizes primitive recursive sequences
of ordinals and the equivalent theorem about primitive recursive Goodstein sequences is expressible in the language of $\PA$ (see Theorem \ref{2.7}).

Barwise in 1977 in the Handbook of Mathematical Logic added an editor's note to the famous paper by Paris and Harrington \cite{ph}:
\begin{quote}Since 1931, the year G\"odel's Incompleteness Theorems were published, mathematicians have been
looking for a strictly mathematical example of an incompleteness in first-order Peano arithmetic, one which is mathematically
simple and interesting and does not require the numerical coding of notions from logic. The first such examples were found
early in 1977, when this Handbook was almost finished.
 \end{quote}
 Barwise describes the problem of finding a natural mathematical incompleteness in Peano arithmetic almost as a ``holy grail problem" of mathematical logic.
 As Goodstein almost found a mathematical example in the 1940s one wonders whether this problem was perceived as so important back then.\footnote{In view of the impact Hilbert's problem list had on mathematics and of how Hilbert's work and ideas furnished the young G\"odel with problem to solve, one might guess that
 if one of the luminaries of mathematical logic had declared the importance of this problem in the 1940s, the young
 ones would have leapt at this chance and followed Gentzen's and Goodstein's lead.}  In his paper Goodstein identifies as his main
 objective to determine which initial cases of Gentzen's restricted
  ordinal theorem can be proven by  finistist means. Here initial cases refer to the ordinals (in Goodstein's notation)
  $\vartheta_n$ ($\vartheta_0=\omega$, $\vartheta_{k+1}=\omega^{\vartheta_k}$) and the pertaining assertions $P(n)$ that all descending sequences below
  $\vartheta_n$ are necessarily finite. Interestingly enough, Goodstein claimed  that
   {\em ``$P(n)$ is capable of a finite constructive proof for any assigned $n$''} (\cite[p.39]{goodstein}).
   Bernays referred to this  claim in his letter from 29 September 1943 \cite{bernays2} when he wrote {\em ``I think, that the methodological difficulties appearing already in the case of $\omega^{\omega^{\omega}}$ , ..., will induce you to speak in a more reserved form about it."} This time Goodstein did not heed Bernays' advice.

In the next section we give an account of Goodman's theorem and illuminate its origins in Gentzen's work. We also give two independence results from $\PA$ that by and large can
be credited to Gentzen (Theorem \ref{2.7}) and Goodstein (Theorem \ref{2.8}), respectively. We leave it to the reader to assay whether they meet
Barwise's criteria of being mathematically
simple and interesting and not requiring the numerical coding of notions from logic.

 In the last section we give an elementary proof of Kirby's and Paris' 1982 result \cite{kp} that the termination of special Goodstein
sequences induced by the shift function is not provable in $\PA$.
Yet another proof was presented by Cichon \cite{cichon} in 1983.
Our
main technique consists in making descending sequences of ordinals $<\varepsilon_0$ slow.
We like to think that this elementary proof could have been found in the 1940's or at least 1950's.
Such a thought could be considered to be unfair to the logicians who, after
a lot of hard technical work, established this independence result. This is not our intention and we like to stress that their techniques were certainly not available before the 1970's. On the other hand, we definitely think that Goodstein and Gentzen deserve at least
some credit for ``their" independence results. Another question that seems to be
relevant in this context is the following: Could it be that the problem
of
finding statements independent from $\PA$ did not occupy center stage in mathematical logic before
the 1970's, thereby accounting for their late arrival?\footnote{For what it's worth, here is some anecdotal evidence. Around 1979, Diana Schmidt
proved that Kruskal's theorem elementarily implies that the ordinal representation system for $\Gamma_0$ is well-founded \cite{schmidt}. She even wrote (p. 61) that she didn't  know of any applications of her result to proof theory. This is quite surprising since in conjunction with proof-theoretic work of Feferman and Sch\"utte from the 1960's it immediately implies the nowadays celebrated result that Kruskal's theorem
is unprovable in predicative mathematics.}

\section{Cantor normal forms}
Let $\varepsilon_0$ be the least ordinal $\beta$ such that $\omega^{\beta}=\beta$.
Every ordinal $0<\alpha<\varepsilon_0$ can be written in a unique way
as \begin{eqnarray}\label{1}\alpha&=&\omega^{\alpha_1}\cdot k_1+\ldots+\omega^{\alpha_n}\cdot k_n\end{eqnarray}
where $\alpha>\alpha_1>\ldots>\alpha_n$ and $0<k_1,\ldots,k_n<\omega$.
This we call the {\bf Cantor normal form} of $\alpha$.  By writing
$\alpha\nf\omega^{\alpha_1}\cdot k_1+\ldots+\omega^{\alpha_n}\cdot k_n$ we shall convey that (\ref{1}) obtains.

The ordinals $\alpha_i$ with $\alpha_i\ne 0$ can also be written in Cantor normal form with yet smaller exponents.
 As this process terminates after finitely many steps every ordinal $<\varepsilon_0$ can be represented in a unique way
 as a term over the alphabet $\omega,+,\cdot,0,1,2,3,\ldots$ which we call its {\bf complete Cantor normal form}.

 In what follows we identify ordinals $<\varepsilon_0$ with their representation in complete Cantor normal form.
 Henceforth, unless indicated otherwise, ordinals are assumed to be smaller than $\varepsilon_0$ and will be denoted
 by lower case Greek letters. By $|\alpha|$ we denote the length of $\alpha$ in complete Cantor normal form  (viewed as a string of symbols). More precisely, if $\alpha\nf\omega^{\alpha_1}\cdot k_1+\ldots+\omega^{\alpha_n}\cdot k_n$
 we define
  $$|\alpha|\,=\,\max\{|{\alpha_1}|,\ldots,|\alpha_n|, k_1,\ldots ,k_n\}+1.$$
  By $\ko(\alpha)$ we denote the highest integer coefficient that appears in $\alpha$, i.e., inductively this can be defined by letting  $\ko(0)=0$ and $$\ko(\alpha)=\max\{\ko(\alpha_1),\ldots,\ko(\alpha_n),k_1,\ldots,k_n\}$$
 where
 $\alpha\nf\omega^{\alpha_1}\cdot k_1+\ldots+\omega^{\alpha_n}\cdot k_n$.

 There is a similar Cantor normal form for positive integers $m$ to any base $b$ with $b\geq 2$,
 namely we can express $m$ uniquely in the form \begin{eqnarray}\label{2} m&=& b^{n_1}\cdot k_1+\ldots +b^{n_r}\cdot k_r\end{eqnarray}
 where $m>n_1>\ldots>n_r\geq 0$ and $0<k_1,\ldots, k_r<b$. As each $n_i>0$ is itself of this form we can repeat
 this procedure, arriving at what is called the {\bf complete $b$-representation} of $m$. In this way
 we get a unique representation of $m$ over the alphabet $0,1,\ldots,b,+,\cdot$.

 For example
$7\,625\,597\,485\,157= 3^{27}\cdot 1+ 3^4\cdot 2 +  3^1\cdot 2+3^0\cdot 2
 =   3^{3^3}+ 3^{3+1}\cdot 2 +  3^1\cdot 2+2$.

 \begin{dfn}{\em Goodstein \cite{goodstein} defined operations mediating between ordinals $<\varepsilon_0$ and
 natural numbers.

 For naturals $m>0$ and $c\geq b\geq 2$ let $\Su b c m $ be the integer resulting from $m$ by replacing the base $b$ in the complete $b$-representation of $m$ everywhere by $c$. For example $\Su 34 {34}=265$, since $34=3^3+3\cdot 2+1$
 and $4^4+4\cdot 2+1=265$.

 For any ordinal $\alpha$ and natural $b\geq 2$ with $b>\ko(\alpha)$ let $\Tu b\alpha$ be the integer resulting from $\alpha$ by replacing
 $\omega$ in the complete Cantor normal form of $\alpha$ everywhere by $b$.
 For example $$\Tu 3{\omega^{\omega+1}+\omega^2\cdot 2+\omega\cdot 2 +1}= 3^{3+1}+3^2\cdot 2+3\cdot 2+1=106.$$

 Conversely, for naturals $m\geq 1$ and $b\geq 2$ let $\Tus bm$ be the ordinal obtained from the complete $b$-representation of $m$ by replacing the base $b$ everywhere with $\omega$.
 Thus $\Tus 3{106}=\omega^{\omega+1}+\omega^2\cdot 2+\omega\cdot 2 +1$
 and $$\Tus 3{34}=\Tus 4{\Su 34{34}}=\Tus 4{265}=\omega^{\omega}+\omega\cdot 2+1.$$

 We also set $\Tus b0=0$ and $\Tu b 0:=0$.
}\end{dfn}

Goodstein's main insight was that given two ordinals $\alpha,\beta<\varepsilon_0$ one could replace the
base $\omega$ in their complete Cantor normal forms by a sufficiently large number $b$ and the resulting
natural numbers $\Tu b\alpha$ and $\Tu b\beta$ would stand in the same ordering as $\alpha$ and $\beta$.
This is simply a consequence of the fact that the criteria for comparing ordinals in Cantor normal form are
the same as for natural numbers in base $b$-representation, as spelled out in the next Lemma.

\begin{lem}\label{vergleich} $\phantom{P}$ \begin{itemize}
\item[(i)] Let $\alpha\nf\omega^{\alpha_1}\cdot k_1+\ldots+\omega^{\alpha_r}\cdot k_r$ and
$\beta\nf\omega^{\beta_1}\cdot l_1+\ldots+\omega^{\beta_s}\cdot k_s$.  Then $\alpha<\beta$
if and only if either of the following obtains:
\begin{enumerate}
\item There exists $0<j\leq \min(r,s)$ such that $\alpha_i=\beta_i$ and $k_i=l_i$ for $i=1,\ldots,j-1$ and
$\alpha_j<\beta_j$, or $\alpha_j=\beta_j$ and $k_j<l_j$.
\item $r<s$ and $\alpha_i=\beta_i$ and $k_i=l_i$ hold for all $1\leq i\leq r$.
\end{enumerate}
\item[(ii)] Let $b\geq 2$, $n=b^{a_1}\cdot k_1+\ldots+b^{a_r}\cdot k_r$ and
$m=b^{a_1'}\cdot l_1+\ldots+b^{a_s'}\cdot l_s$ be $b$-representations of integers $n$ and $m$, respectively.
  Then $n<m$
if and only if either of the following obtains:
\begin{enumerate}
\item There exists $0<j\leq \min(r,s)$ such that $a_i=a_i'$ and $k_i=l_i$ for $i=1,\ldots,j-1$ and
$a_j<a_j'$, or $a_j=a_j'$ and $k_j<l_j$.
\item $r<s$ and $a_i=a_i'$ and $k_i=l_i$ hold for all $1\leq i\leq r$.
\end{enumerate}
 \end{itemize}
 \end{lem}

\begin{lem}\label{lem1} Let $m,n,b$ be naturals, $b\geq 2$, and  $\alpha,\beta$ be ordinals with $\ko(\alpha),\ko(\beta)<b$.
\begin{itemize} \item[(i)] $\Tu b{\Tus bm}=m$.
\item[(ii)] $\Tus b{\Tu b\alpha}=\alpha$.
\item[(iii)] $\alpha<\beta \Leftrightarrow \Tu b\alpha<\Tu b \beta$.
\item[(iv)] $m<n \Leftrightarrow \Tus bm<\Tus bn$.
\end{itemize}
\end{lem}
\prf (i) and (ii) are obvious. (iii) and (iv) follow from Lemma \ref{vergleich}. \qed

\begin{dfn}\label{good} {\em Given any natural number  $m$ and non-decreasing function $$f:{\mathbb N}\to {\mathbb N}$$ with  $f(0)\geq 2$
 define
            $$ m_0^f =m,\;\;\ldots\; \;,
            m_{i+1}^f= \Su {f(i)}{f(i+1)}{m^f_i}\minpunkt 1$$
            where  $k\minpunkt 1$ is the predecessor of $k$ if $k>0$, and $k\minpunkt 1=0$ if $k=0$.

            We shall call $(m_i^f)_{i\in {\mathbb N}}$ a {\bf Goodstein sequence}.
            Note that a sequence $(m_i^f)_{i\in {\mathbb N}}$ is uniquely determined by $f$ once we fix
            its starting point $m=m_0^f$.

            The case when $f$ is just a shift function has received special attention.
            Given any $m$ we define $m_0=m$ and $m_{i+1}:=\Su {i+2}{i+3}{m_i}\minpunkt 1$
          and call $(m_i)_{i\in {\mathbb N}}$ a {\bf special Goodstein sequence}.
          Thus $(m_i)_{i\in {\mathbb N}}= (m^{\mathrm{id}_2}_i)_{i\in {\mathbb N}}$, where  $\mathrm{id}_2(x)=x+2$.
          Special Goodstein sequences can differ only with respect to their starting points.
          They give rise to a recursive function $f_{good}$
          defined as follows: $f_{good}(m)$ is the least $i$ such that $m_i=0$ where
          $(m_i)_{i\in\N}$ is the special Goodstein sequence starting with $m_0=m$.
}
\end{dfn}

\begin{thm}[Goodstein 1944]\label{thm1} Every Goodstein sequence terminates, i.e. there exists
$k$ such that $m_i^f=0$ for all $i\geq k$.
\end{thm}
\prf If $m^f_i\ne 0$ one has $$\Tus {f(i)} {m^f_i}=\Tus {f(i+1)}{\Su{f(i)}{f(i+1)}{m^f_i}}>\Tus {f(i+1)}{m^f_{i+1}}$$
by Lemma \ref{lem1}(iv) since  $\Su{f(i)}{f(i+1)}{m^f_i}=m^f_{i+1}+1$. Hence, as there are no infinitely descending ordinal sequences, there must exist a $k$ such that $m^f_k=0$. \qed

The statement of the previous theorem is not formalizable in $\PA$. However, the corresponding statement about termination of special Goodstein sequences is expressible in the language of $\PA$ as a $\Pi_2$ statement.
It was shown to be unprovable in $\PA$ by Kirby and Paris in 1982 \cite{kp} using model-theoretic tools.
\cite{kp} prompted  Cichon  \cite{cichon} to find a different (short) proof that harked back to older proof-theoretic work of Kreisel's \cite{kreisel} from 1952 which identified the so-called $<\varepsilon_0$-recursive functions as the provably recursive functions
of $\PA$. Other results pivotal to \cite{cichon} were ordinal-recursion-theoretic classifications of Schwichtenberg \cite{schwichtenberg} and Wainer \cite{wainer} from around 1970 which showed that the latter class of recursive functions consists exactly of those elementary in one of the fast growing functions $F_{\alpha}$ with $\alpha<\varepsilon_0$.
As $F_{\varepsilon_0}$ eventually dominates any of these functions it is not provably
total in $\PA$. Cichon verified that  $F_{\varepsilon_0}$ is elementary in the function $f_{good}$ of Definition \ref{good}. Thus termination of special Goodstein sequences is not provable in $\PA$.

Returning to Goodstein, he established a connection between sequences of natural numbers and descending sequences of ordinals.
Inspection of his proof shows that, using the standard scale of reverse mathematics,  it can be carried out in the weakest system, $\mathbf{RCA}_0$, based on of recursive comprehension (see \cite{SOSA}).

\begin{thm}[Goodstein 1944] Over $\mathbf{RCA}_0$ the following are equivalent:
\begin{itemize} \item[(i)] Every Goodstein sequence terminates.
\item[(ii)] There are no infinitely descending sequences of ordinals $$\varepsilon_0>\alpha_0>\alpha_1>\alpha_2>\ldots.$$
    \end{itemize}
    Of course, when we speak about ordinals  $<\varepsilon_0$ in $\mathbf{RCA}_0$ we mean Cantor normal forms.
\end{thm}
\prf ``(ii)$\Rightarrow$(i)" follows from Theorem \ref{thm1}.
For the converse, assume (i) and, aiming at a contradiction, suppose we have a strictly descending sequence of ordinals $\varepsilon_0>\alpha_0>\alpha_1>\alpha_2>\ldots$.
 Define a function $f:\mathbb N\to \mathbb N$ by letting $f(i)=\max\{\ko(\alpha_0),\ldots,\ko(\alpha_i)\}+1$.
 $f$ is non-decreasing. Let $m:= \Tu {f(0)}{\alpha_0}$. We claim that
    \begin{eqnarray}\label{t2a} \Tu {f(i)}{\alpha_i}&\leq& m^f_i\end{eqnarray}  and to this end use induction on $i$. It's true for $i=0$ by definition.
    Inductively assume $m^f_i\geq \Tu{f(i)}{\alpha_i}$. Then $$\Su {f(i)}{f(i+1)}{m^f_i}\geq \Su {f(i)}{f(i+1)}{\Tu{f(i)}{\alpha_i}}$$ and hence
    \begin{eqnarray}\label{t2b}&& \Su {f(i)}{f(i+1)}{m^f_i}\geq \Su {f(i)}{f(i+1)}{\Tu{f(i)}{\alpha_i}}=\Tu{f(i+1)}{\alpha_i}>\Tu{f(i+1)}{\alpha_{i+1}}
    \end{eqnarray} where the last inequality holds by Lemma \ref{lem1}(iii) since
     $\alpha_{i+1}<\alpha_i$ and $$\ko(\alpha_{i+1}),\ko(\alpha_i)<f(i+1).$$
     From (\ref{t2b}) we conclude that $m^f_{i+1}= \Su {f(i)}{f(i+1)}{m^f_i}\minpunkt 1\geq \Tu{f(i+1)}{\alpha_{i+1}}$,
     furnishing the induction step.

     Since $m^f_k=0$ for a sufficiently large $k$, (\ref{t2a}) yields that $\Tu{f(k)}{\alpha_k}=0$, and hence
     $\alpha_k=0$, contradicting $\alpha_k>\alpha_{k+1}$. \qed
  Whereas it's not possible to speak about arbitrary Goodstein sequences in $\PA$,
  one can certainly formalize the notion of a primitive recursive sequence of naturals in this theory. As a result of the proof of the previous Theorem we have:

  \begin{cor}\label{prim-goodstein} Over $\mathbf{PA}$ the following are equivalent:
\begin{itemize} \item[(i)] Every primitive recursive Goodstein sequence terminates.
\item[(ii)] There are no infinitely descending primitive recursive sequences of ordinals $$\varepsilon_0>\alpha_0>\alpha_1>\alpha_2>\ldots.$$
    \end{itemize}
   \end{cor}

A very coarse description of Gentzen's result \cite{gentzen} one often finds is that he showed that transfinite induction up to $\varepsilon_0$ suffices to prove the consistency
of first order number theory (also known as Peano arithmetic, $\PA$).
What Gentzen actually did is much more subtle. He defined a reduction procedure on derivations (proofs) and showed that if
successive application of a reduction step on a given derivation always leads to a non-reducible derivation in finitely many steps, then the consistency of $\PA$ follows. The latter he ensured by assigning ordinals to derivations in such a way that a reduction step applied to a reducible derivation results in a derivation with a smaller ordinal.
Let us explain in more detail how this is done in the later \cite{gentzen38} which uses the sequent calculus.
Firstly, he defined an assignment $\ord$ of ordinals to derivations of $\PA$ such for every derivation $D$ of $\PA$ in his sequent calculus, $\ord(D)$ is an ordinal $<\varepsilon_0$. He then defined a reduction procedure
$\RR$ such that whenever $D$ is a derivation of the empty sequent in $\PA$ then $\RR(D)$ is another derivation of the empty sequent in $\PA$ but with a smaller ordinal assigned to it, i.e., \begin{eqnarray}\label{gentz}&&\ord(\RR(D))<\ord(D).\end{eqnarray}
Moreover,
both $\ord$ and $\RR$ are primitive recursive functions and only finitist means are used in showing  (\ref{gentz}).
As a result, if $\PRWO(\varepsilon_0)$ is the statement that there are no infinitely descending primitive recursive
sequences of ordinals below $\varepsilon_0$, then the following are immediate consequences of Gentzen's work.

\begin{thm}[Gentzen 1936, 1938]\label{2.7} $\phantom{NN}$ \begin{itemize}
\item[(i)] The theory of primitive recursive arithmetic, $\PRA$, proves that
$\PRWO(\varepsilon_0)$ implies the consistency of $\PA$.
\item[(ii)] Assuming that $\PA$ is consistent, $\PA$ does not prove $\PRWO(\varepsilon_0)$.
\end{itemize}
\end{thm}
\prf For (ii), of course, one  invokes G\"odel's second incompleteness theorem. \qed
So it appears that an attentive reader could have
inferred the following from \cite{gentzen,gentzen38,goodstein} in 1944:

\begin{thm}\label{2.8} Termination of primitive recursive Goodstein sequences is not
provable in $\PA$.
\end{thm}
\prf Use Theorem \ref{2.7} (ii) and Corollary \ref{prim-goodstein}.\qed
\section{Slowing down}
The key to establishing that already termination of special Goodstein sequences is beyond $\PA$ is
to draw on Theorem \ref{2.7} and to show that infinite descending sequences can be made slow.
This technology was used in a paper by Simpson \cite[Lemma 3.6]{simpson-kruskal} where it is credited to Harvey Friedman. It would be good to know where this padding technique was used for the first time.
\begin{dfn} {\em Addition of ordinals $\alpha+\beta$ is usually defined by transfinite recursion on $\beta$.   For ordinals given in complete Cantor normal form
addition can be defined explicitly. We set $\alpha+0:=\alpha$ and $0+\alpha:=\alpha$.
Now let $\alpha,\beta$ be non-zero ordinals, where
 $\alpha\nf\omega^{\alpha_1}\cdot k_1+\ldots+\omega^{\alpha_r}\cdot k_r$ and
$\beta\nf\omega^{\beta_1}\cdot l_1+\ldots+\omega^{\beta_s}\cdot l_s$.
If $\alpha_1<\beta_1$ then $\alpha+\beta:=\beta$.
Otherwise there is a largest $1\leq i\leq r$ such that $\alpha_{i}\geq\beta_1$.
If $\alpha_{i}=\beta_1$, then $$\alpha+\beta:=\omega^{\alpha_1}\cdot k_1+\ldots+\omega^{\alpha_{i-1}}\cdot k_{i-1}+\omega^{\beta_1}\cdot (k_{i}+l_1)
+\omega^{\beta_2}\cdot l_2+\ldots+\omega^{\beta_{l_s}}\cdot l_s.$$
If $\alpha_{i}<\beta_1$, then $$\alpha+\beta:=\omega^{\alpha_1}\cdot k_1+\ldots+\omega^{\alpha_{i}}\cdot k_{i}+\omega^{\beta_1}\cdot l_1
+\ldots+\omega^{\beta_{l_s}}\cdot l_s.$$
With the help of addition we can also explicitly define multiplication
$\omega^{\alpha}\cdot \beta$ as follows:
$\omega^{\alpha}\cdot 0:=0$. If $\beta\nf\omega^{\beta_1}\cdot l_1+\ldots+\omega^{\beta_s}\cdot l_s$ then
$$\omega^{\alpha}\cdot \beta:= \omega^{\alpha+\beta_1}\cdot l_1+\ldots+\omega^{\alpha+\beta_s}\cdot l_s.$$
We shall use $\omega\cdot\beta$ to stand for $\omega^1\cdot\beta$.
 }\end{dfn}
 Next we recall an elementary result that was known in the 1950's (e.g., \cite{grze}).
 \begin{lem}\label{grz1} For a function $\ell:\N \to \N$ define
  $\ell^0(l)=l$ and $\ell^{k+1}(l)=\ell(\ell^k(l))$.
  The {\em Grzegorczyk hierarchy} $(f_l)_{l\in\mathbb N}$ is generated by the functions $f_0(n)=n+1$ and $f_{l+1}(n)= (f_l)^n(n)$.

   For every primitive recursive function $h$ of arity $r$
  there is an $n$ such $h(\vec{x}\,)\leq f_n(\max(2,\vec x\,))$ holds for all
  $\vec x=x_1,\ldots,x_r$.
  \end{lem}
  \prf The proof proceeds by induction on the generation of primitive recursive
  function, using properties of the hierarchy $(f_l)_{l\in\mathbb N}$.
  It is straightforward but a bit tedious. We shall give it in the appendix.
  \qed

 \begin{lem}[$\PA$]\label{3.2} Let $f:{\mathbb N}\to {\mathbb N}$ be primitive recursive. Then there exists a primitive recursive function $g:{\mathbb N}^2\to \omega^{\omega}$ such that
 \begin{enumerate}
 \item[(1)] $g(n,m)>g(n,m+1)$ whenever $m<f(n)$.
 \item[(2)] There exists a constant $K$ such $|g(n,m)|\leq K\cdot (n+m+1)$ holds for all $n,m$.
 \end{enumerate}
 \end{lem}
 \prf By Lemma \ref{grz1} is suffices to show this for any $f=f_l$ in the
 in the hierarchy $(f_l)_{l\in\mathbb N}$.
  We will actually obtain a $0<k<\omega$ such that $g:{\mathbb N}^2\to \omega^{k}$. To find $g$ we proceed by induction on $l$.

\noindent
{\bf Base Case}: $f(n)=n+1$.
Define $g$ by $$g(n,m)=(n+2)\minpunkt  m.$$
{\bf Induction Step}: Let $g:{\mathbb N}^2\to \omega^{k}$ satisfy the conditions (1) and (2) for $f$, and
let $f'$ be defined by diagonalizing over $f$, i.e.,  $f'(k)=f^k(k)$, where $f^0(l)=l$ and $f^{k+1}(l)=f(f^k(l))$.
If $m<f'(n)$ define $g'(n,m)$ by letting
$$g'(n,m)=\omega^k\cdot (n-i)+g(f^i(n),j),$$
where $i$ and $j$ are the unique integers such that
$$m=f(n)+f^2(n)+\ldots +f^i(n)+j,$$
$i<n$ and  $j<f^{i+1}(n)$. If $m\geq f'(n)$ set $g'(n,m)=0$.

We first show that $g'$ satisfies requirement (1)
for $f'$. So suppose $m<f'(n)$. Let
$m=f(n)+f^2(n)+\ldots +f^i(n)+j$ with $j<f^{i+1}(n)$.
We distinguish two cases. If also $j+1<f^{i+1}(n)$, 
then $$g'(n,m+1)=\omega^k\cdot(n-i)+g(f^i(n),j+1)<\omega^k\cdot(n-i)+g(f^i(n),j)=
g'(n,m)$$ holds by the inductive assumption on $g$ and $f$ since $j<f(f^i(n))$.
 The other possible case is that $j+1=f^{i+1}(n)$ and then we have
 $$g'(n,m+1)=\omega^k\cdot(n-(i+1))+g(f^{i+1}(n),0)<\omega^k\cdot(n-i)+g(f^i(n),j)=
g'(n,m)$$ since $\omega^k\cdot(n-(i+1))+\omega\leq\omega^k\cdot(n-i)$ as $k>0$.

$g'$ also satisfies requirement (2)  for $f'$ since
\begin{eqnarray*}|g'(n,m)|&\leq& \mathrm{constant}\cdot  n+\mathrm{constant} \cdot (f^i(n)+m+1)\\
                           &\leq & \mathrm{constant}\cdot (n+m+1).\end{eqnarray*}
 \qed

 \begin{cor}[$\PA$]\label{3.3} From a given primitive recursive strictly descending sequence $\varepsilon_0>\beta_0>\beta_1>\beta_2>\ldots$ one can construct a slow primitive recursive strictly descending sequence
 $\varepsilon_0>\alpha_0>\alpha_1>\alpha_2>\ldots$, where slow means that
  there is a constant $K$ such that $$|\alpha_i|\leq K\cdot(i+1)$$ holds for all $i$.
  \end{cor}
  \prf By the previous Lemma let $g:{\mathbb N}^2\to \omega^{\omega}$ be chosen such that
 $g(n,m)>g(n,m+1)$ for every $m<|\beta_{n+1}|$ and $|g(n,m)|\leq K\cdot (n+m+1)$ holds for all $n,m$.
Now set
$$\alpha_j\,=\,\omega^{\omega}\cdot \beta_n+g(n,m)$$
where $j=|\beta_0|+|\beta_1|+\ldots+|\beta_{n}|+m$ for $m<|\beta_{n+1}|$.
For such $j$ one computes that
 \begin{eqnarray*}|\alpha_j|&\leq& \mathrm{constant}\cdot |\beta_n|+\mathrm{constant}\cdot (n+m+1)\\
                               &\leq& \mathrm{constant}\cdot(j+1).\end{eqnarray*}
 We also need to determine $\alpha_i$ for $i<|\beta_0|$. For instance let
$$\alpha_i\,=\,\omega^{\omega}\cdot \beta_0+|\beta_0|+1-i$$
for $i<|\beta_0|$. Clearly we can choose a constant $K_0$ such that $|\alpha_i|\leq K_0\cdot(i+1)$ for $i<|\beta_0|$.                              \qed

\begin{thm}[$\PA$]\label{3.4} Let $\varepsilon_0>\alpha_0>\alpha_1>\alpha_2>\ldots$ be a slow primitive recursive descending sequence of ordinals, i.e.,
 there is a constant $K$ such that $|\alpha_i|\leq K\cdot(i+1)$ for all $i$.
 Then there exists  a primitive recursive descending sequence $\varepsilon_0>\beta_0>\beta_1>\beta_2>\ldots$
 such that $\ko(\beta_r)\leq r+1$ for all $r$.
 \end{thm}
 \prf Obviously $K>0$. Let $\omega_0=1$ and $\omega_{n+1}=\omega^{\omega_n}$. As $\alpha_0<\varepsilon_0$ we find $s<\omega$ such that
 $\omega\cdot \alpha_0<\omega_s$ and $K<s$.
  Now put 
  $$\beta_j\,:=\,\sum_{i=0}^{K-1-j}\omega_{s-i}$$
  for $j=0,\ldots,K-1$, and
   $$\beta_{K\cdot(n+1)+i}\,:=\,\omega\cdot\alpha_n+(K-i)$$ for $n<\omega$ and $0\leq i<K$.
   By construction, $\beta_r>\beta_{r+1}$ for all $r$.
   As $\ko(\omega_r)=1$ for all $r$, one has $\ko(\beta_j)=1$ for all $j=0,\ldots,K-1$.
    Moreover, as $\ko(\alpha_n)\leq |\alpha_n|\leq K\cdot(n+1)$, it follows that
     $$\ko(\beta_{K\cdot(n+1)+i})=\ko(\omega\cdot\alpha_n+(K-i))\leq K\cdot(n+1)+1,$$
     since multiplying by $\omega$ increases the coefficients by at most 1.
    As a result, $\ko(\beta_r)\leq r+1$ for all $r$. \qed

    \begin{lem}[$\PA$]\label{3.5} Let $\varepsilon_0>\beta_0>\beta_1>\beta_2>\ldots$ be a primitive recursive descending sequence of ordinals
    such that $\ko(\beta_n)\leq n+1$.
    Then the special Goodstein sequence $(m_i)_{i\in\mathbb N}$ with $m_0=\Tu 2{\beta_0}$ and
    $m_{i+1}=\Su{i+2}{i+3}{m_i}\minpunkt  1$ does not terminate.
    \end{lem}
    \prf 
     We claim that \begin{eqnarray}\label{a}&& m_k\,\geq\,\Tu {k+2}{\beta_k}\end{eqnarray}
     holds for all $k$.

     For $k=0$ this holds by definition. Assume this to be true for $i$, i.e., $m_i\geq \Tu {i+2}{\beta_i}$.
     Let $\delta=\Tus {i+2}{m_i}$. Since $\ko(\beta_i)<i+2$ it follows from Lemma \ref{lem1}(iii)
     that $\delta\geq \beta_i$, and hence $\delta>\beta_{i+1}$. As $\ko(\delta),\ko(\beta_{i+1})<i+3$ it follows from Lemma \ref{lem1}(iii) that $\Tu{i+3}{\delta}>\Tu{i+3}{\beta_{i+1}}$. Thus, since $$m_{i+1}=\Su{i+2}{i+3}{m_i}\minpunkt  1=\Tu{i+3}{\delta}\minpunkt  1,$$ we arrive at
      $m_{i+1}\geq \Tu{i+3}{\beta_{i+1}}$ as desired.

       (\ref{a}) entails $m_k\ne 0$ for all $k$.  \qed

       In sum, what we have done amounts to an elementary proof of the following result due to Kirby and Paris
       \cite[Theorem 1(ii)]{kp}:

       \begin{cor}\label{3.6} The statement that any special Goodstein sequence terminates is not provable in $\PA$.
       \end{cor}
       \prf Let $GS$ be the statement that every special Goodstein sequence terminates. Arguing in $\PA$ and assuming $GS$, we obtain from Lemma \ref{3.5}, Theorem 3.4 and Corollary \ref{3.3} that there is no infinite primitive recursive descending sequence
       of ordinals below $\varepsilon_0$, i.e., $\PRWO(\varepsilon_0)$. However, by Theorem \ref{2.7} the latter is
       not provable in $\PA$. \qed

\section{Appendix}
It remains to prove Lemma \ref{grz1}. To this end the following is useful.

 \begin{lem}\label{grz1a} Recall that for a function $h:\N \to \N$ we defined
  $h^0(l)=l$ and $h^{k+1}(l)=h(h^k(l))$.
  Also recall that the hierarchy $(f_l)_{l\in\mathbb N}$ is generated by the functions $f_0(n)=n+1$ and $f_{l+1}(n)= (f_l)^n(n)$.
  We shall write $f^n_l$ rather than $(f_l)^n$.

  Let $f$ be any of the functions $f_l$
  in this hierarchy. Then $f$ satisfies the following properties:
  \begin{itemize}
  \item[(i)] $f(x)\geq x+1$ if $x>0$.
  \item[(ii)] $f^z(x)\geq x$ for all $x,z$.
  \item[(iii)] If $x<y$ then $f(x)<f(y)$ and $f^z(x)<f^z(y)$.
  \item[(iv)] $f_{l+1}(x)\geq f_l(x)$ whenever $x>0$.
  \end{itemize}
  \end{lem}
  \prf (i),(ii),(iii) will be proved simultaneously by induction on $l$.
  (i) and (iii) are obvious for $f=f_0$ and (ii) follows via a trivial induction on $z$. Now assume that (i),(ii),(iii) hold for $f_k$ and $l=k+1$.
  For $x>0$ one then computes $$f_l(x)=f_k^x(x)=f_k(f^{x-1}_k(x))\geq f_k(x)\geq x+1$$
  using the properties for $f_k$.  (ii) follows from this by induction on $z$.
  As to (iii), note that
  $$f_l(x+1)=f_k^{x+1}(x+1)=f_k(f_k^x(x+1))>f_k(f_k^x(x))\geq f^x_k(x)=f_l(x),$$
  using the properties for $f_k$, and thus
  (iii) follows by straightforward inductions on $y$ and $z$.

  If $x>0$, then $f_{l+1}(x)=f^x_l(x)=f_l(f^{x-1}_l(x))\geq f_l(x)$ by (ii) and (iii).
  \qed

  {\bf Proof of Lemma \ref{grz1}}: We want to prove
that for or every primitive recursive function $h$ of arity $r$
  there is an $n$ such $h(\vec{x}\,)\leq f_n(\max(2,\vec x\,))$ holds for all
  $\vec x=x_1,\ldots,x_r$.

  We show this by induction on the generation of the primitive recursive functions.
  Clearly for all $n$ we have $h(\vec{x})\leq f_n(\max(2,\vec x\,))$ by Lemma \ref{grz1a}(i)
  if $h$ is any of the initial functions $x\mapsto 0$, $\vec{x}\mapsto x_i$, and
  $x\mapsto x+1$.

  Now let $h$ be defined by $h(\vec x\,)=g(\varphi_1(\vec x\,),\ldots,\varphi_s(\vec x\,))$ and assume that the assertion holds for $g,\varphi_1,\ldots,\varphi_s$.
  By Lemma \ref{grz1}(iv) we can then pick an $n$ such that
  $g(\vec y\,)\leq f_n(\max(2,\vec y\,))$ and $\varphi_i(\vec x\,)\leq f_n(\max(2,\vec y\,))$ hold for
  all $\vec y,\vec x$ and $1\leq i\leq s$. As a result,
  \begin{eqnarray*} h(\vec x\,) &\leq &f_n(\max(2,f_n(\max(2,\vec x\,))))=f_n(f_n(\max(2,\vec x\,)))= \\ &&
  f_n^2(\max(2,\vec x\,))\leq f_n^{\max(2,\vec x\,)}(\max(2,\vec x\,))=f_{n+1}
  (\max(2,\vec x\,)),\end{eqnarray*}
  showing that $f_{n+1}$ is a majorant for $h$.

  Now suppose $h$ is defined by primitive recursion from $g$ and $\varphi$ via
   $h(\vec x,0)=g(\vec x\,)$ and $h(\vec x,y+1)=\varphi(\vec x,y,h(\vec x,y))$
   and that $f_n$ majorizes $g$ and $\varphi$, i.e., $g(\vec x\,)\leq f_n(\max(2,\vec x\,))$ and $\varphi(\vec x,y,z)\leq f_n(\max(2,\vec x\,))$.
   We claim that \begin{eqnarray}\label{rec} h(\vec x,y) &\leq & f_n(\max(2,\vec x,y)).
   \end{eqnarray}
   We prove this by induction on $y$. For $y=0$ we have $h(\vec x,y)=g(\vec x\,)\leq
   f_n(\max(2,\vec x\,))=f^1_n(\max(2,\vec x\,))$. For the induction step we
   compute \begin{eqnarray*} h(\vec x,y+1) &=& \varphi(\vec x,y,h(\vec x,y))\leq
   f_n(\max(2,\vec x,y,h(\vec x,y)))\\
   &\leq& f_n(\max(2,\vec x,y,f_n^{y+1}(\max(2,\vec x,y))))= f_n(f_n^{y+1}(\max
   (2,\vec x,y))) \\
   &=& f_n^{y+2}(\max(2,\vec x,y))\end{eqnarray*} where the
   second ``$\leq$" uses the inductive assumption and the penultimate ``$=$" uses
   Lemma \ref{grz1a}.

   From the claim (\ref{rec}) we get with Lemma \ref{grz1a}, letting $w=\max(2,\vec x,y)$, that
   \begin{eqnarray*} h(\vec x,y)&\leq& f_n^{y+1}(\max(2,\vec x,y))\leq f_n^{w+1}(w)=
   f_n(f_n^w(w))=f_n(f_{n+1}(w))\\
   &\leq& f_{n+1}(f_{n+1}(w))=f^2_{n+1}(w)\leq f_{n+1}^w(w)=f_{n+2}(w).\end{eqnarray*}
   As a result, $h(\vec x,y)\leq f_{n+2}(\max(2,\vec x,y))$. \qed

\paragraph{Acknowledgement:}
I thank Jan von
Plato for having shared his archival findings about the
Goodstein-Bernays correspondence with me even before his own work on it
\cite{plato} has been  published.

I also acknowledge
 support by the EPSRC of the UK through
 Grant No. EP/G029520/1.


\begin{thebibliography}{99}
\bibitem{hb} J. Barwise (ed.): {\em Handbook of Mathematical Logic},
(North-Holland, Amsterdam, 1977).
\bibitem{bernays} P. Bernays: {\em Letter to Goodstein, dated September 1st, 1942}, Bernays collection of the ETH Z\"urich.
\bibitem{bernays2} P. Bernays: {\em Letter to Goodstein, dated September 29th, 1943}, Bernays collection of the ETH Z\"urich.
 \bibitem{cichon} A. Cichon:
 {\em A short proof of two recently discovered independence results using recursion theoretic methods}, Proc. Amer. Math. Soc. 87 (1983) 704--706.
 \bibitem{gentzen}G. Gentzen: {\em Die Widerspruchsfreiheit der reinen
Zahlentheorie.} Mathematische Annalen 112 (1936) 493--565.
\bibitem{gentzen38} G. Gentzen: {\em Neue Fassung des Widerspruchsfreiheitsbeweises f\"ur die reine Zahlentheorie}, Forschungen zur Logik und zur Grundlegung der exacten Wissenschaften, Neue Folge 4 (Hirzel, Leipzig, 1938) 19--44.
\bibitem{goodstein} R.L. Goodstein 1944:  {\em On the restricted ordinal theorem}, Journal of Symbolic Logic 9 (1944)
33--41.

\bibitem{grze} A. Grzegorczyk: {\em Some classes of recursive functions}. Rozprawy
Mate No. IV (Warsaw, 1953).
 \bibitem{kp} L. Kirby, J. Paris: {\em Accessible independence results for Peano arithmetic}, Bull. London Math. Soc. 14 (1982) 285--293.
\bibitem{kreisel} G. Kreisel: {\em On the interpretation of non-finitist proofs II.} Journal of Symbolic Logic 17 (1952) 43--58.
\bibitem{ph} J. Paris, L. Harrington: {\em A mathematical incompleteness in Peano arithmetic}. In:
J. Barwise (ed.): {\em Handbook of Mathematical Logic},
(North-Holland, Amsterdam, 1977) 1133--1142.
\bibitem{plato} Jan von Plato: {\em G\"odel, Gentzen, Goodstein: The magic sound of a G-string}. To appear in the {\em Mathematical Intelligencer}.
    \bibitem{schmidt} D. Schmidt: {\em Well-partial orderings and their maximal order types.} (Habilitationsschrift, Universit\"at Heidelberg, 1979).
\bibitem{sc}K. Sch\"utte: {\em Proof Theory.} Springer 1977.

\bibitem{schwichtenberg} H. Schwichtenberg: {\em Eine Klassifikation der $\varepsilon_0$-rekursiven Funktionen}. Zeitschrift f\"ur mathematische Logik und Grundlagen der Mathematik 17 (1971) 61--74.

\bibitem{simpson-kruskal} S.G. Simpson: {\em Nichtbeweisbarkeit von gewissen kombinatorischen Eigenschaften endlicher B\"aume}, Archiv f\"ur mathematische Logik 25 (1985) 45--65.

\bibitem{SOSA} S.G. Simpson: {\em Subsystems of Second Order Arithmetic}, second edition,
(Cambridge University Press, 2009).

\bibitem{wainer} S.S. Wainer: {\em A classification of the ordinal recursive functions.} Archiv f\"ur Mathematische Logik und Grundlagenforschung
13 (1970) 136--153.

\end{thebibliography}
\end{document}